\documentclass{amsart}
\usepackage{graphicx}
\usepackage{amssymb}
\usepackage{graphicx}
\usepackage{amsfonts}
\usepackage{amsmath}
\usepackage{amsthm}
\usepackage{esint}
\usepackage{mathabx}
\usepackage{mathtools}
\usepackage{nccmath}
\usepackage[colorlinks]{hyperref}

\newtheorem{thm}{Theorem}[section]
\newtheorem{cor}[thm]{Corollary}
\newtheorem{lem}[thm]{Lemma}
\newtheorem{prop}[thm]{Proposition}
\theoremstyle{definition}

\theoremstyle{remark}
\newtheorem{rem}[thm]{Remark}
\theoremstyle{remark}
\newtheorem{ex}[thm]{Example}

\newcommand{\ch}{\mathop{{\rm ch}}}
\renewcommand{\chi}{{ 1}}
\renewcommand{\H}{{ \mathcal H}}

\newcommand{\spt}{\mathop{\rm spt}}

\newcommand{\p}{\partial}
\newcommand{\R}{\mathbf R}

\begin{document}

\title[]{ Discontinuities in optimal transport between nonconvex domains}
\author{Robert McCann}
\address{Departments of Mathematics and of Economics, University of Toronto}
\email{mccann@math.toronto.edu}

\author{Arghya Rakshit}
\address{Department of Mathematics, University of Toronto}
\email{\tt arghya.rakshit@utoronto.ca}



\thanks{
$^*$Robert McCann's work was supported in part by the Canada Research Chairs program CRC-2020-00289 and Natural Sciences and Engineering Research Council of Canada Discovery Grant RGPIN--2026-04906. 
Figure \ref{fig:Vregion} was generated with ChatGPT. All mathematical ideas and results are the authors' own.
    \copyright \today}

 
\begin{abstract}
We develop a new approach to studying discontinuity sets of  measure-preserving convex gradient maps between nonconvex domains. Such maps optimize the Euclidean distance squared. Our idea is to introduce a new `phantom' optimal transport problem, that transports the missing mass from the convex hull of the target to the `excess' Monge--Ampère measure of the minimal convex extension of the Brenier potential. We show  that the closure of the  discontinuity set of the original map  contains the support of this phantom target. Moreover, inside a planar source domain these two sets actually coincide;   on $n$-dimensional domains they coincide $\mathcal{H}^{n-1}$-a.e. The Hausdorff exponent $n-1$ is sharp in every dimension $n\geq3$. For planar polygonal targets, we establish a quantitative global length bound for the closure of the discontinuity set depending only on the target geometry, valid for bounded positive densities.
\end{abstract}
\maketitle




\section{Introduction}\label{Intro}

 Among volume-preserving maps between given domains in $\R^n$,  the Euclidean distance squared is optimized by a convex gradient \cite{Brenier91}, \cite{McCann95}. Caffarelli showed convexity of the support of the target implies smoothness of the map in the source interior \cite{C1}. Discontinuities may arise if the target is not convex (see  e.g. Caffarelli \cite{C1}, Figalli \cite{F1}, Kitagawa--McCann \cite{KM1}, Chen--Liu \cite{ChenLiu25}, Chen--Li--Liu \cite{CLL1}, Jhaveri \cite{Jhaveri19}, Mooney--Rakshit \cite{MR1} and Rakshit--Sen~\cite{RS1}).  By discontinuities of the optimal transport map $\nabla u$, we mean the set of points where  its potential $u$  is not differentiable. Our goal is to develop a method to detect discontinuities of optimal transport maps between nonconvex domains, especially in the plane.\\

Let us fix some notation. For a set $A\subset \R^n$, $\ch(A)$ denotes the convex hull of $A$. For a convex function $f$, we use $\partial f$ to denote the subgradient of $f$,  and $\nabla f$ to denote its gradient. Thus $\partial f(x) = \{\nabla f(x)\}$ holds $\H^n$-a.e.,  where $\H^k$ is the $k$-dimensional Hausdorff measure on $\R^n$. In particular, $\H^n$ is Lebesgue measure  or volume, which we also denote by $|A| := \H^n(A)$.\\

 Our setting is as follows. Let \(U,V\subset\R^n\) be bounded, open sets, which can be taken to have unit volume without loss of generality \cite{GangboMcCann00}. Brenier's theorem \cite{Brenier91} provides a globally defined convex potential \(u:\R^n\to\R\)  whose gradient gives a volume-preserving map between $U$ and $V$:
\[
(\nabla u)_\#(\chi_U\,dx)=\chi_V\,dy.
\]
We define the minimal convex extension of \(u|_U\) by
\[
\bar u(x):=\sup_{\substack{z\in U\\p\in\partial u(z)}}
\{u(z)+p\cdot(x-z)\}.
\]
Then \(\bar u=u\) on \(U\), and \(\bar u\leq w\) for every convex \(w:\R^n\to\R\) satisfying \(w=u\) on \(U\). Notice that we do not assume convexity of the domain \(U\).

Let \(v:=\bar u^*\) be the Legendre transform of \(\bar u\). By definition, $v$ is convex.  It is real valued throughout $\ch(V)$, but identically $+\infty$ outside $\overline{\ch(V)}$ due to the minimality of $\bar u$. We can use it to define an optimal transport map $\nabla v$ that takes uniform density on $\ch(V)\setminus V$  to some measure which is a priori unknown. 

We refer to $\ch(V)\setminus V$ as the cavity in the target. Our first theorem identifies this unknown measure in terms of the Monge-Amp\`ere measure $\mathcal M \bar u$. Recall for each Borel set $A \subset \R^n$, the Monge-Amp\`ere measure of $u$ is defined by
$$
(\mathcal M \bar u)(A) = |\p \bar u(A)|.
$$

The following theorem is an immediate consequence of \cite[Lemma 4.1]{McCann97}, as is the generalization to nonuniform probability densities discussed in Remark \ref{R:nonuniform1} below.

\begin{thm}[Convex hull of the target accounts for excess Monge-Amp\`ere mass]
\label{phantomtheorem}
Let $U,V \subset \R^n$ be bounded, open, unit volume sets. Then 
\[
 (\nabla v)_\# (1_{\ch(V)} \, d\H^n) = d\mathcal M \bar u.
\]
Hence
\begin{equation}\label{uniform}
     (\nabla v)_{\#}\left(\chi_{\ch(V)\setminus V}\,d\H^n\right)
=     d\mathcal M\bar u-\chi_U\,d\H^n. 
\end{equation}
\end{thm}

Let's denote the phantom source measure above by $d\nu=\chi_{\ch(V)\setminus V}\,d\mathcal H^n$ and  the phantom target measure  by $d\mu=d\mathcal{M}\bar{u}-\chi_U\,d\mathcal H^n$. We call  $\nabla v$ the \emph{phantom optimal transport map}. The goal of this paper is to relate the discontinuity set of the optimal transport map $\nabla u$ to this  phantom optimal transport problem.

Let $S_{\bar{u}}$ be the discontinuity set of $\nabla\bar{u}$, or equivalently,  the set of points in $\R^n$ where $\bar{u}$ is not differentiable. Define $\Sigma$ to be the closure of  $S_{\bar{u}}$. We define the support of \(\mu\) by
\[
\operatorname{spt}\mu
:=
\left\{
x\in \R^{n}:\mu(B_r(x))>0 \text{ for every } r>0
\right\}.
\]
 The transport identity above provides the starting point for our analysis. Our main result identifies the support of the excess Monge--Amp\`ere measure with the closed discontinuity set inside a planar source domain. In higher dimensions, this identification holds up to an \(\H^{n-1}\)-null set:

\begin{thm}[Relating map discontinuities to phantom Monge-Amp\`ere mass]
\label{main}
Under the hypotheses and notation of Theorem~\ref{phantomtheorem}, $\mathcal{H}^{n}(\partial V)=0$  implies
\[
\spt\mu\subset\Sigma
\qquad\text{and}\qquad
\H^{n-1}\bigl((\Sigma\cap U)\setminus\spt\mu\bigr)=0.
\]
If, in addition, \(n=2\), then
\[
\Sigma\cap U=\spt\mu\cap U.
\]
\end{thm}
{The higher-dimensional conclusion reflects a distinction between discontinuities of the transport map and concentration of excess Monge--Amp\`ere mass. Gradient discontinuities away from the support of singular Monge--Amp\`ere mass were constructed in  Mooney \cite{MooneyY} with one of us \cite{MRdirac} for equations with point masses. Our phantom transport framework gives this phenomenon a geometric interpretation for transport between uniform densities: the point masses arise from the collapse of target cavities, while discontinuities also occur away from the support of the collapsed mass. These examples therefore reveal an obstruction to extending the planar identification of the closed discontinuity set with the support of phantom mass to higher dimensions. 
 Building on work of Caffarelli \cite{CaffarelliMcCann10} and Kitagawa \cite{KM1} with McCann, Chen--Liu \cite{ChenLiu25} establish regularity of the codimension-one interface arising when the target consists of two convex domains separated by a hyperplane. Our framework does not require such a decomposition of the target, and our examples exhibit discontinuities of higher codimension that do not separate the source into two transport phases. Our sharpness construction quantifies the extent of this failure and shows that the Hausdorff exponent \(n-1\) in Theorem \ref{main} is optimal (see Section \ref{sharp}). The closed discontinuity set inside \(U\) can still be characterized exactly by considering multiple contacts with the target boundary (see Subsection \ref{sec:exactcharacterization}).}

\begin{rem}
 Since \(\partial\bar u(U\setminus S_{\bar u})\subset\overline V\), applying \eqref{eq:phantom-subgradient} and using \(\mathcal H^n(\partial V)=0\), we obtain
\begin{equation}\label{eq:no-mass-on-closure-points}
\mu(U\setminus S_{\bar u})=0,
\qquad
\mu\bigl((\Sigma\setminus S_{\bar u})\cap U\bigr)=0.
\end{equation}
Thus, in dimension two, the additional closure points of the discontinuity set inside \(U\) carry no phantom mass, although they belong to \(\operatorname{spt}\mu\). How large this set can be is a separate question, answered for polygonal targets in Theorem~\ref{partial}.
\end{rem}

\noindent Next, we discuss some generalization of the theorems above:
\begin{rem}[Nonuniform probability densities]
\label{R:nonuniform1}
    Theorem \ref{phantomtheorem} is valid for source and target densities $f$ and $g$ supported on $U$ and $V$ respectively, and 
\begin{equation}\label{probability densities}
    \int_U f=\int_V g=1 \quad\text{with}\quad f>0\ {\rm and}\ g >0
\end{equation}
     a.e. in $U$ and $V$, provided the right-hand side of \eqref{uniform} is replaced by
    \[
    dM \bar u
    -\frac{f(x)}{g(\nabla u(x))}\chi_U(x)\,dx.
    \]
 This is because
    \[
\begin{aligned}
\int_V \varphi(\nabla v(y))\,dy
&=
\int_V \frac{\varphi(\nabla v(y))}{g(y)}\,g(y)\,dy \\
&=
\int_U \varphi(x)\frac{f(x)}{g(\nabla u(x))}\,dx.
\end{aligned}
\]
Therefore,
\[
(\nabla v)_{\#}\bigl(\mathbf 1_V\,dy\bigr)
=
\frac{f(x)}{g(\nabla u(x))}\mathbf 1_U(x)\,dx.
\]
\end{rem}

\begin{rem}[Uniform bounds]
If the source and target densities $f$ and $g$ satisfy \eqref{probability densities} and
\begin{equation}\label{uniform bounds}
0<\lambda\leq f,g\leq\Lambda<\infty,
\end{equation}
then the conclusion of Theorem~\ref{main} remains valid, with
\[
d\mu \;=\; dM_{\bar u}-\frac{f(x)}{g(\nabla u(x))}\,\chi_U\,dx
\]
defined as in the previous remark. The boundedness of the densities away from $0$, and $\infty$ is needed to apply Caffarelli's theory for the Brenier potential away from $\Sigma\cap U$(see \cite{C1}).
\end{rem}

\begin{rem}[Discontinuity versus singular set in dimension $2$]
 If we are working with smooth source and target densities that are bounded away from $0$ and $\infty$, then we will have smoothness of the optimal transport map {in the interior of the domain except on} the closure of the discontinuity sets. In this case the singular (non-smooth) set of the potential and the discontinuity set of the optimal transport map coincide in $U$.
\end{rem}

As a geometric application of Theorem~\ref{main}, we show that, under a separation assumption on the target cavities, the closed discontinuity set inside a planar source decomposes into separated portions corresponding to those cavities (see Proposition \ref{prop:separation}).\\

Next, we study some interesting examples. First, we construct a case where $\Sigma$ is a line that goes to infinity in both directions. In Example \ref{ex:explicit}, we show that the support of the phantom optimal transport map can lie  outside $U$ or even outside $\overline{U}$. This has particular relevance to the semigeostrophic theory of the atmosphere and oceans, which was devised by Hoskins \cite{Hoskins75} as a model for pressure front formation and evolution;  see also Cullen \cite{Cullen06}  and the references therein.  In this context,  the target $V$ represents a snapshot of an advected potential vorticity patch relative to a fixed fluid domain $U$,  and the discontinuities $\Sigma = \overline {S_{\bar u}}$ represent a snapshot of the corresponding pressure fronts.  If $V$ evolves smoothly from a convex to a nonconvex shape,  the pressure front forms outside $U$ initially (by the regularity shown for small nonconvex perturbations $V$ of convex targets of Chen and Figalli \cite{CF1}), and is expected to later enter the physical domain from its boundary,  if the evolving target $V$ becomes sufficiently non-convex.

Finally,  for $n=2$ we study the discontinuity set of the optimal transport map when the target is a polygon. In dimension two, Figalli
 showed the discontinuous set inside $U$ has Hausdorff dimension at most $1$ and described in it further detail in \cite{F1},  establishing smoothness of $u$  in the interior of $U \setminus \Sigma$. Chen, Li and Liu proved that if the target is a polygon, then  
$ U \cap \Sigma$ is locally a smooth one-dimensional curve away from finitely many points (see \cite{CLL1}). We build on these developments to prove the following sharp partial regularity result when the target is a polygon. In contrast to the existing local descriptions of the discontinuity set, we obtain a global quantitative estimate for the length of its closure:

\begin{thm}[Polygonal targets yield finite length discontinuities]\label{partial}
Let $U\subset\R^2$ be a bounded open domain, and $V$ the interior of a polygon in $\R^2$. Let $f$ and $g$ be densities on $U$ and $V$, respectively, with equal mass and
\[
0<\lambda\le f,g\le\Lambda<\infty
\]
almost everywhere. For $\Sigma$ as above, we have
\[
\mathcal{H}^1(\Sigma\cap U)
\le C_V\,\operatorname{diam}(U),
\]
 where $C_V$ depends only on the geometry of $V$, independently of the densities and their bounds.
\end{thm}
In fact, $C_V$ admits an explicit upper bound determined solely by the locations of the vertices of $V$ (see Remark~\ref{rem:explicit}).\\

{Theorem~\ref{main} identifies the closed discontinuity set itself inside \(U\) as the support of the phantom measure, while Theorem~\ref{partial} controls its total length. For polygonal targets and uniform densities, we further describe how this mass is distributed inside \(U\): the phantom measure has an integrable density with respect to arclength along the regular curves, together with finitely many atoms (see Corollary~\ref{cor:structure}). This gives a finer description than an estimate of the size of the singular set alone, by also describing the mass it carries. Even a set of zero \(\H^1\)-measure can carry positive phantom mass, so estimates of the size of the singular set alone do not describe how this mass is distributed. Each atom at \(x\) has mass \(|\partial\bar u(x)|\). Thus a single singular point may carry positive phantom mass, even though it has zero \(\H^1\)-measure.}\\

\noindent The paper is organized as follows. In Section \ref{proofs} we first prove Theorem \ref{phantomtheorem}. Then we state some critical lemmas, and use them to justify the proof of Theorem \ref{main}. We prove a separation result for the singularities generated by distinct target cavities  in Section \ref{sec:separation}. In Section~\ref{examples}, we discuss all the examples mentioned above. In Section~\ref{partialsec}, we prove Theorem \ref{partial}. Finally, in Section \ref{future}, we mention some open questions for further study.

\section{Proof of Theorem \ref{phantomtheorem} and \ref{main}}\label{proofs}

\subsection{Proof of Theorem \ref{phantomtheorem}}

\begin{proof} 
With the identifications $\psi:=v$, $\psi^*=\bar{u}$, $\rho := 1_{\ch(V)}=g$,  and $M:=A$ the first claim becomes a special case of McCann \cite[Lemma 4.1]{McCann97}. On the other hand, since $v=u^*$ in $V$, the map $\nabla v$ takes uniform density on $V$ to $\chi_U\,dx$.  Thus \[ (\nabla v)_{\#}\left(\chi_{\ch(V)\setminus V}\,d\H^n\right) = {d}\mathcal M\bar u-\chi_U\,d\H^n. \] 
\end{proof}

\subsection{Some Results Towards Proof of Theorem \ref{main} }
The first observation is that  the discontinuity set $\Sigma$ accounts for the entire phantom volume contained in the cavity $\ch(V)\setminus V$.

\begin{lem}[Image of continuity points have negligible intersection with cavities]
\label{nomassoutofsing}
If \(B\subset \R^{n}\setminus \Sigma\) is a ball, then
\[
\left|\partial \bar u(B)\cap \big(\ch(V)\setminus V\big)\right|=0.
\]
\end{lem}

\begin{proof} 
Since \(B\cap \Sigma=\emptyset\), the function \(\bar u\) is  differentiable in \(B\). Therefore,
\[ \partial \bar u(B)=\nabla \bar u(B). \] 
 We address the interior, complement and boundary of $\overline U$ in turn. On \(B\cap U\), we have \(\bar u=u\), and the original transport map satisfies \[ (\nabla u)_{\#}(\chi_U\,dx)=\chi_V\,dy. \] Thus, \(\nabla \bar u(B\cap U)\subset \overline V\), up to a set of Lebesgue measure zero. Hence,
\[ \left|\nabla \bar u(B\cap U)\cap \big(\ch(V)\setminus V\big)\right|=0. \]
By the minimality of the convex extension, for every $x\in B\setminus \overline U$ there exists a point $z\in\overline U$, with $z\neq x$, such that the same affine function supports $\bar u$ at both $x$ and $z$. Let $p=\nabla\bar u(x)$. By Legendre duality,
\[
x,z\in\partial\bar u^*(p).
\]
Thus $\bar u^*$ is not differentiable at $p$. Since a convex function is differentiable almost everywhere,
\[
\left|\nabla\bar u(B\setminus\overline U)\right|=0.
\]
For $x\in B\cap\partial U$, the continuity of $\nabla\bar u$ in $B$ and the fact that $\nabla u$ takes values in $V$ almost everywhere in $U$ imply that $\nabla\bar u(x)\in\overline V$.  Thus $\mathcal H^{n}(\p V)=0$ implies
\[
\left|
\nabla\bar u(B\cap\partial U)
\cap\bigl(\ch(V)\setminus V\bigr)
\right|=0.
\]
Consequently,
\[
\left|
\nabla\bar u(B\setminus U)
\cap\bigl(\ch(V)\setminus V\bigr)
\right|=0.
\]
Combining the above estimates gives \[ \left|\partial \bar u(B)\cap \big(\ch(V)\setminus V\big)\right|=0. \]
\end{proof}

\noindent This immediately gives an inclusion of the support of the phantom measure inside the closure of the discontinuity set  which holds in all dimensions:

\begin{prop}[Discontinuities cover excess Monge-Amp\`ere mass]\label{P:Sigmaisbigger}
With \(\mu\) as above, we have
\[
\operatorname{spt}\mu\subset \Sigma.
\]
\end{prop}

\begin{proof}
Let \(x\notin \Sigma\). Since \(\Sigma\) is closed, there exists \(r>0\) such that
\[
B_r(x)\subset \R^{n}\setminus \Sigma.
\]
By Theorem \ref{phantomtheorem}
\[
\mu(A)
=
\bigl|\{y\in \ch(V)\setminus V:\nabla v(y)\in A\}\bigr|.
\]
Since $v=\bar u^*$ is differentiable almost everywhere, Legendre duality gives
\[
\nabla v(y)=x
\quad\Longleftrightarrow\quad
y\in\partial\bar u(x)
\]
for almost every $y$. Hence, up to a null set,
\[
\{y:\nabla v(y)\in A\}=\partial\bar u(A),
\]
and therefore
\begin{equation}\label{eq:phantom-subgradient}
\mu(A)
=
\bigl|\partial\bar u(A)\cap(\ch(V)\setminus V)\bigr|.
\end{equation}
 Thus
\[
\mu(B_r(x))=\bigl|\partial\bar u(B_r(x))\cap(\ch(V)\setminus V)\bigr|=0
\]
and hence \(x\notin \operatorname{spt}\mu\). Therefore
\[
\operatorname{spt}\mu\subset \Sigma.
\]
\end{proof}

We use the previous proposition to prove the main result in this paper:\\
\subsection{Proof for Theorem \ref{main}}

\begin{proof}
First we prove 
\[
\Sigma\cap U=\spt\mu\cap U.
\]
in dimension $2$. Let \(x_0\in U\setminus\spt\mu\). Since \(x_0\notin\spt\mu\) and \(U\) is open, there exists \(r>0\) such that \(B_{2r}(x_0)\Subset U\) and \(\mu(B_{2r}(x_0))=0\). Thus, we have \(\mu(A)=0\) for every Borel set \(A\subset B_{2r}(x_0)\). Since \(d\mu=d\mathcal M\bar u-\chi_U\,d\H^2\), we obtain
\[
(\mathcal M\bar u)(A)=|A|
\qquad\text{for every Borel set }A\subset B_{2r}(x_0).
\]
Thus \(\bar u\) is an Alexandrov solution of \(\det D^2\bar u=1\) in \(B_{2r}(x_0)\). Since \(n=2\), \cite[Lemma 2.3 and Remark 2.4]{Mooney15} implies that \(\bar u\) is strictly convex there. Caffarelli's interior regularity theorem \cite{C1} then gives \(\bar u\in C^{1,\alpha}(B_r(x_0))\) for some \(\alpha\in(0,1)\).

 We conclude that \(S_{\bar u}\cap B_r(x_0)=\emptyset\), and therefore \(x_0\notin\overline{S_{\bar u}}=\Sigma\). Hence \(U\setminus\spt\mu\subset U\setminus\Sigma\), which proves \(\Sigma\cap U\subset\spt\mu\cap U\). Combining this with Proposition~\ref{P:Sigmaisbigger} yields
\[
U\cap\Sigma=U\cap\spt\mu.
\]
Next we prove, without assuming \(n=2\), one still has
\[
\H^{n-1}\bigl((\Sigma\cap U)\setminus\spt\mu\bigr)=0.
\]
Indeed, cover \(U\setminus\spt\mu\) by countably many balls \(B_j\Subset U\setminus\spt\mu\). On each \(B_j\), we have \(d\mathcal M\bar u=d\H^n\), so \(\bar u\) is an Alexandrov solution of \(\det D^2\bar u=1\). By \cite[Theorem~1.1]{Mooney15}, \(\bar u\) is strictly convex in \(B_j\) away from a set \(E_j\) with \(\H^{n-1}(E_j)=0\). At every strictly convex point, Caffarelli's interior regularity theory \cite{C1} gives \(C^{1,\alpha}\)-regularity in a neighborhood. Hence \((\Sigma\cap B_j)\setminus\spt\mu\subset E_j\), and the conclusion follows by taking the countable union.
\end{proof}

\begin{rem}
The proof of the reverse inclusion uses only the measure identity proved in  Theorem~\ref{phantomtheorem}, strict convexity of Alexandrov solutions in dimension two, and Caffarelli's interior regularity theory. It does not  require the connectedness of \(V\),  any boundary regularity of \(V\)  beyond $\mathcal H^{n}(\p V)=0$, or the assumption that \(\operatorname{ch}(V)\setminus V\) has finitely many connected components. 
\end{rem}
\begin{rem}
The only genuinely dimension-dependent step in the proof is the automatic strict convexity of Alexandrov solutions in dimension two. Caffarelli's interior \(C^{1,\alpha}\) regularity theorem applies in every dimension once strict convexity is known. Since \(d\mathcal M\bar u=d\H^n\) on the open set \(U\setminus\spt\mu\), combining \cite[Theorem~1.1]{Mooney15} with Caffarelli's theorem there gives
\[
\H^{n-1}\bigl((\Sigma\cap U)\setminus\spt\mu\bigr)=0 .
\]
Thus the equality in Theorem~\ref{main} holds exactly in dimension two and, in every dimension, up to a set of zero \(\H^{n-1}\)-measure.
\end{rem}

\subsection{Sharpness of Theorem \ref{main}}\label{sharp}
{The planar identity \(\Sigma\cap U=\spt\mu\cap U\) can fail in every dimension \(n\ge3\), even for uniform densities and a ball $U=B$ as source. Indeed, \cite{MRdirac} and \cite{MooneyY}, give a convex function \(u\) satisfying
\[
\mathcal Mu= \mathcal{H}^n+\delta_{e_1}+\delta_{-e_1},
\]
where \(e_1=(1,0,\ldots,0)\in \R^n\). We fix the full asymptotic quadratic polynomial \(Q\) of this construction and require \(u(x)-Q(x)\to0\) as \(|x|\to\infty\), which determines \(u\) uniquely. This convex function \(u\) is smooth outside \(L:=[-e_1,e_1]\) and nondifferentiable at every point of \(L\). Choose a ball \(U\) containing \(L\) compactly and set \(V:=\nabla u(U\setminus L)\). Then \(V\) is bounded and open, \(\H^n(\partial V)=0\), and \(\nabla u\) transports \(\chi_U\,dx\) to \(\chi_V\,dy\). The minimal convex extension \(\bar u\) of \(u|_U\) agrees with \(u\) in \(U\), so
\[
\mu|_U=\delta_{e_1}+\delta_{-e_1},
\qquad
\spt\mu\cap U=\{-e_1,e_1\}
\subsetneq [-e_1,e_1]=\Sigma\cap U.
\]
The two Dirac masses arise from collapsing positive-volume portions of the target cavity to the endpoints. Indeed, the phantom transport identity gives
\[
\left|\left\{y\in\ch(V)\setminus V:
\nabla v(y)=\pm e_1\right\}\right|=1.
\]
Thus each endpoint receives one unit of cavity mass. In contrast, every \(x\in(-e_1,e_1)\) has a neighborhood \(B_r(x)\Subset U\setminus\{\pm e_1\}\) with \(\mu(B_r(x))=0\), although \(u\) is not differentiable at \(x\). The cavity collapse therefore accounts for the endpoint atoms, but does not detect the discontinuities along the open segment, which persist without any phantom mass nearby. }\\

We next show that the exponent \(n-1\) in
Theorem~\ref{main} cannot be improved.
Mooney \cite[Section~4]{Mooney15} constructed solutions of \(\mathcal M_u=\H^n\) in $B_1$ that fail to be strictly convex on sets of Hausdorff dimension arbitrarily close to \(n-1\). However, failure of strict convexity does not imply nondifferentiability. We modify his comparison barriers to force two distinct supporting slopes along a Cantor family of affine segments, producing genuine discontinuities of the gradient while preserving \(\mathcal M_u=\H^n\) in $B_1$ .

The resulting construction gives, for every \(n\geq3\) and \(0<\eta<1\), transport between uniform densities on a ball \(U\) and a bounded open target \(V\), with \(\H^n(\partial V)=0\), such that
\[
\spt\mu\cap U=\varnothing,
\qquad
\H^{n-1-\eta}(\Sigma\cap U)=\infty.
\]
Thus the discontinuity set can have Hausdorff dimension arbitrarily close to \(n-1\) without any phantom mass inside the source. Related constructions of discontinuity sets are
discussed in Rakshit and Sen \cite{RS1}. We first establish the required nondifferentiability estimate.

\begin{lem}[Sharpness of dimension bound]\label{lem:large-nondifferentiability-set}
For every \(n\ge3\) and \(0<\eta<1\), there exists a convex, locally Lipschitz function \(u:B_1\subset\mathbb R^n\to\mathbb R\) such that
\[
\mathcal M_u=\H^n\quad\text{in }B_1,
\qquad
\H^{n-1-\eta}
\bigl(\{x\in B_{3/4}:u\text{ is not differentiable at }x\}\bigr)
=\infty.
\]
\end{lem}
\begin{proof}
Choose the parameter \(\delta>0\) in
\cite[Section~4]{Mooney15} sufficiently small that
\[
\frac32\delta<\eta,
\qquad
\alpha:=2-\delta,\qquad
\beta:=1+\frac{\delta}{4-3\delta}.
\]
Thus \(1<\alpha,\beta<2\) and
\(\alpha^{-1}+\beta^{-1}=3/2\).
Retain a compact smaller copy \(S\Subset(-1/2,1/2)\) of the Cantor set constructed there.
Then
\[
\dim_{\mathcal H}S=1-\frac32\delta>1-\eta,
\qquad
\H^{1-\eta}(S)=\infty.
\]
Let \(h\) be the convex function denoted by \(v\) in that construction.
For each \(a\in S\), choose the supporting line
\[
\ell_a(s):=h(a)+p_a(s-a),\qquad p_a\in\partial h(a),
\]
provided by the  construction. By \cite[Section~4]{Mooney15}, there are constants
\(c_*>0\) and \(r_0\in(0,1/2)\), independent of \(a\in S\), such that
\[
h(s)-\ell_a(s)\geq c_*|s-a|^\alpha
\qquad\text{whenever }|s-a|\leq r_0.
\]
Since \(h-\ell_a\) is convex, nonnegative, and vanishes at \(a\),
for \(|s|\leq1\) and \(r_0<|s-a|\leq3/2\) we have
\[
h(s)-\ell_a(s)
\geq c_*r_0^{\alpha-1}|s-a|
\geq c_*\left(\frac{2r_0}{3}\right)^{\alpha-1}|s-a|^\alpha.
\]
Thus, setting \(c_0:=c_*(2r_0/3)^{\alpha-1}>0\), we obtain
\[
h(s)-\ell_a(s)\geq c_0|s-a|^\alpha
\qquad(a\in S,\ |s|\leq1).
\]
Set \(\Omega:=D_1\times(-1,1)\) to be a three dimensional cylinder, that is $D_1:=\{x_1^2+x_2^2<1\}$.
Using the explicit barrier in \cite[Section~4]{Mooney15},
choose constants \(A>0\) and \(\kappa>0\) so that
\[
w(y):=A\bigl(|y_1|^\alpha+|y_2|^\beta\bigr)
             (1+\kappa y_3^2)
\]
is convex and satisfies
\[
\mathcal M_w\geq\H^3
\qquad\text{in }D_2\times(-1,1).
\]
This follows by first taking \(\kappa\) sufficiently small and
then \(A\) sufficiently large. Moreover, \(w\geq0\),
\(w(0,0,t)=0\), and
\[
w(y)\leq C_1\bigl(|y_1|^\alpha+|y_2|^\beta\bigr)
\qquad
\text{on }\overline{D_2}\times[-1,1],
\]
where \(C_1:=A(1+\kappa)\).\\
For \(a\in S\), define
\[
\phi_a(x):=w(x_1-a,x_2,x_3).
\]
Since \(|a|<1/2\), the translated cylinder lies in
\(D_2\times(-1,1)\). Hence
\[
\mathcal M_{\phi_a}\geq\H^3\quad\text{in }\Omega,
\qquad
\phi_a(a,0,t)=0\quad(|t|<1).
\]
Moreover, \(|x_2|^\beta\leq|x_2|\) on \(\overline\Omega\), so
\[
\begin{aligned}
\phi_a(x)
&\leq C_1\bigl(|x_1-a|^\alpha+|x_2|^\beta\bigr)\\
&\leq C_0\bigl(h(x_1)-\ell_a(x_1)+|x_2|\bigr)
\qquad(x\in\overline\Omega),
\end{aligned}
\]
where \(C_0:=C_1\max\{c_0^{-1},1\}\) is independent of \(a\).\\
Fix \(c>0\), choose \(C\ge C_0+c\), and use \cite[Lemma~4.1]{Mooney15} to obtain the convex solution \(u_3\in C(\overline\Omega)\) of
\[
\mathcal M_{u_3}=\H^3\quad\text{in }\Omega,\qquad
u_3=C\bigl(h(x_1)+|x_2|\bigr)\quad\text{on }\partial\Omega.
\]
For \(a\in S\), replace the lower comparison barrier by
\[
b_a(x):=C\ell_a(x_1)+\phi_a(x)+c|x_2|.
\]
Then \(b_a\) is convex and \(\mathcal M_{b_a}=\mathcal M_{\phi_a}\) on each of the half-spaces \(x_2>0\) and \(x_2<0\), since the added terms are affine there. As \(\{x_2=0\}\) has zero volume and \(\mathcal M_{b_a}\) is nonnegative, it follows that \(\mathcal M_{b_a}\ge\H^3\) in \(\Omega\). On \(\partial\Omega\), we have
\[
C\bigl(h(x_1)+|x_2|\bigr)-b_a(x)
\ge (C-C_0)(h-\ell_a)(x_1)
   +(C-C_0-c)|x_2|
\ge0.
\]
Comparison therefore gives \(u_3(x)\ge C\ell_a(x_1)+c|x_2|\). Since \(u_3(a,0,\pm1)=Ch(a)\), convexity also gives \(u_3(a,0,t)\le Ch(a)\) for \(|t|<1\). Consequently,
\[
u_3(a,0,t)=Ch(a),\qquad
(Cp_a,c,0),\ (Cp_a,-c,0)\in\partial u_3(a,0,t).
\]
Since \(c>0\), these supporting slopes are distinct. Thus \(u_3\) is nondifferentiable throughout \(S\times\{0\}\times(-1,1)\).

For \(n=3\), take \(u=u_3|_{B_1}\). For \(n>3\), define
\[
u(x,z):=u_3(x)+\frac12|z|^2,\qquad
(x,z)\in B_1\subset\mathbb R^3\times\mathbb R^{n-3}.
\]
The product formula for the subdifferential gives \(\mathcal M_u=\H^n\). In either case, \(u\) is convex and finite, hence locally Lipschitz, and is nondifferentiable on \(S\times\{0\}\times(-r,r)^{n-2}\). For \(r>0\) sufficiently small, this product is contained in \(B_{3/4}\). Since \(\H^{1-\eta}(S)=\infty\), the slicing inequality for Hausdorff measure implies that the product has infinite \(\H^{n-1-\eta}\)-measure.

\end{proof}
We now apply the lemma to get an optimal transport problem. 
\begin{thm}\label{thm:sharp}
For every \(n\geq3\) and \(0<\eta<1\), there exist bounded open sets \(U,V\subset\mathbb R^n\) of equal positive volume, with \(U\) a ball and \(\mathcal H^n(\partial V)=0\), such that the optimal transport map \(T=\nabla\bar u\) transporting \(\chi_U\,dx\) to \(\chi_V\,dy\) has associated phantom measure \(\mu\) and closed discontinuity set \(\Sigma\) satisfying
\[
\operatorname{spt}\mu\cap U=\emptyset,
\qquad \text{and}\qquad
\mathcal H^{n-1-\eta}(\Sigma\cap U)=\infty.
\]
\end{thm}

\begin{proof}
Fix \(n\ge3\) and \(0<\eta<1\), take \(u\) as above, and set \(U=B_{3/4}\). Let \(G\subset B_1\) be the open set where \(u\) is locally smooth, and define \(V:=\nabla u(U\cap G)\). By \cite[Theorem~1.1]{Mooney15} and interior regularity, \(|B_1\setminus G|=0\). Convexity and \(D^2u>0\) on \(G\) make \(\nabla u:U\cap G\to V\) a diffeomorphism with Jacobian one. Thus \(V\) is bounded and open, and \((\nabla u)_\#(1_U\,dx)=1_V\,dy\).\\
Closedness of the subdifferential graph gives
\[
\partial V\subset
\partial u\bigl(\partial U\cup(U\setminus G)\bigr).
\]
Since \(\partial U\cup(U\setminus G)\) has zero volume and \(\mathcal M_u=\H^n\), it follows that
\[
\H^n(\partial V)
\le\mathcal M_u\bigl(\partial U\cup(U\setminus G)\bigr)=0.
\]
The minimal convex extension \(\bar u\) agrees with \(u\) in \(U\), so the associated phantom measure satisfies \(\mu|_U=0\). Hence \[\spt\mu\cap U=\varnothing .\]
Lemma~\ref{lem:large-nondifferentiability-set} therefore gives
\[
\H^{n-1-\eta}
\bigl((\Sigma\cap U)\setminus\spt\mu\bigr)=\infty.
\]
A common dilation of the source and target normalizes their volumes to one. Since \(\eta>0\) can be arbitrarily small, no smaller Hausdorff exponent can replace \(n-1\) in Theorem~\ref{main}, even for uniform densities and a ball as source. This gives us the following sharpness result for Theorem \ref{main}:

\end{proof}

\subsection{Exact Characterization of the Singular Set}\label{sec:exactcharacterization}
{The phantom measure quantifies the excess Monge–Ampère mass arising from the collapse of target cavities. As we saw above, its support need not capture the entire closed discontinuity set in higher dimensions. We next give an exact characterization of this set in terms of multiple contacts with the target boundary. The Proposition below is a consequence of \cite{C1, FigalliKim10}.
\begin{prop}\label{prop:boundary-characterization}
Assume the hypotheses and notation of
Theorem~\ref{phantomtheorem}. In every dimension,
\begin{equation}\label{eq:boundarycontacts}
S_{\bar u}\cap U
=
\left\{x\in U:
\partial\bar u(x)\cap\partial V
\text{ contains at least two points}\right\}.
\end{equation}
Consequently,
\begin{equation}\label{eq:closurebounarycontacts}
\Sigma\cap U
=
U\cap\overline{
\bigcup_{\substack{p,q\in\partial V\\p\neq q}}
\bigl(\partial v(p)\cap\partial v(q)\bigr)
}.
\end{equation}
\end{prop}
\begin{proof}
For \(x\in S_{\bar u}\cap U\), Caffarelli's extreme-point argument in \cite{C1} places every extreme point of \(\partial\bar u(x)\) in \(\overline V\). By \cite[Proposition 3.4]{FigalliKim10}, if
\(\partial\bar u(x)\cap V\neq\emptyset\), then \(\bar u\) is differentiable at \(x\), a contradiction.
Thus every extreme point of \(\partial\bar u(x)\) lies on \(\partial V\). The non-singleton compact convex set \(\partial\bar u(x)\) therefore has at least two distinct extreme points on \(\partial V\). Conversely, two distinct boundary contacts
imply nondifferentiability. This proves
\eqref{eq:boundarycontacts}.

By Legendre duality,
\[
p,q\in\partial\bar u(x)
\quad\Longleftrightarrow\quad
x\in\partial v(p)\cap\partial v(q).
\]
Taking closures relative to the open set \(U\) proves \eqref{eq:closurebounarycontacts}.
\end{proof}
}
{ 
\begin{rem}
The boundary characterization above locates the discontinuities, but does not quantify the excess Monge--Amp\`ere mass concentrated on them. The phantom transport framework in Theorem \ref{main} provides this finer information through the measure \(\mu\), which records how the collapsed cavity mass is distributed over the closed discontinuity set.
\end{rem}}

\section{Target cavities do not interact in the plane}\label{sec:separation}

For planar domains, we also observe that the discontinuities containing the excess Monge-Amp\`ere mass corresponding to different target cavities remain disjoint:

\begin{prop}[Target cavities do not interact in the plane]
\label{prop:separation}
Under the hypotheses and notation of Theorem~\ref{main}, assume \(n=2\). Let \(\ch(V)\setminus\overline V=T_1\sqcup T_2\) with $\overline T_1 \cap \overline T_2 =\emptyset$ and define
\[
\mu_i:=(\nabla v)_\#(\chi_{T_i}\,d\H^2),
\qquad
\Sigma_i:=\spt\mu_i\cap U,
\qquad i=1,2.
\]
Then
\[
\Sigma\cap U=\Sigma_1\sqcup\Sigma_2.
\]
Moreover, \(\operatorname{dist}(\Sigma_1,\Sigma_2)>0\) whenever both sets are nonempty.
\end{prop}

\begin{proof}
We have \(d\mathcal Mv\geq\chi_V\,dy\). Suppose \(x\in\spt\mu_1\cap\spt\mu_2\). By the definition of support and the pushforward identity, we can choose differentiability points \(p_k\in T_1\) and \(q_k\in T_2\) such that
\[
\lim_{k\to\infty} \nabla v(p_k) =x =
\lim_{k\to\infty} \nabla v(q_k) 
\]
After passing to subsequences,
\[
p_k\to p\in\overline{T_1},
\qquad
q_k\to q\in\overline{T_2}.
\]
Since \(v=\bar u^*\) is a closed convex function, closedness of its subdifferential graph gives
\[
x\in\partial v(p)\cap\partial v(q).
\]
Thus \(v\) agrees on \([p,q]\) with a supporting affine function, which we denote by \(\ell\).

Since \(\overline{T_1}\cap\overline{T_2}=\emptyset\), the connected segment \([p,q]\) cannot be contained in \(\overline{T_1}\cup\overline{T_2}\). Hence there exists \(a\in(p,q)\setminus(\overline{T_1}\cup\overline{T_2})\). By convexity of \(\ch(V)\), for sufficiently small \(r>0\) there is an open half-ball \(B_r^+(a)\) with diameter \(I_r\subset(p,q)\) such that
\[
B_r^+(a)\subset
\ch(V)\setminus(\overline{T_1}\cup\overline{T_2}).
\]
Since \(\ch(V)\setminus\overline V=T_1\sqcup T_2\) and \(\H^2(\partial V)=0\), the set \(V\) has full measure in \(B_r^+(a)\). Therefore \(d\mathcal Mv\geq dy\) in \(B_r^+(a)\).

Extend \(v|_{B_r^+(a)}\) to \(B_r(a)\) by setting it equal to \(\ell\) on \(B_r(a)\setminus B_r^+(a)\). This extension is convex since \(v\geq\ell\) and \(v=\ell\) on \(I_r\). It is affine on \(I_r\), while its Monge--Amp\`ere measure is bounded below by Lebesgue measure on \(B_r^+(a)\). After translation, rotation, and rescaling, this contradicts \cite[Lemma~6.1]{MR1}. Hence
\[
\spt\mu_1\cap\spt\mu_2=\emptyset.
\]

Finally, since \(\H^2(\partial V)=0\),
\[
\chi_{\ch(V)\setminus V}
=\chi_{T_1}+\chi_{T_2}
\qquad\text{a.e.},
\]
and therefore
\[
\mu=\mu_1+\mu_2,
\qquad
\spt\mu=\spt\mu_1\cup\spt\mu_2.
\]
Theorem~\ref{main} yields
\[
\Sigma\cap U
=\spt\mu\cap U
=\Sigma_1\sqcup\Sigma_2.
\]
Since \(U\) is bounded and \(\overline{\Sigma_i}\subset\spt\mu_i\), the two closures are compact and disjoint. Their distance, and hence \(\operatorname{dist}(\Sigma_1,\Sigma_2)\), is positive whenever both sets are nonempty.
\end{proof}

\section{Some Explicit Examples}\label{examples}
The next example illustrates why the identity in Theorem \ref{main} must be restricted to $U$:
\begin{ex}\label{caffarelli}
Consider Caffarelli's classical example \cite{C1}. Take \(U=B_1\subset\R^2\) and
\[
V=(B_1^+ +e_1)\cup (B_1^- -e_1),
\qquad \text{where}\;
B_1^\pm:=B_1\cap\{\pm x>0\}.
\]
Up to a measure zero set, $\ch(V)\setminus V$ is the rectangle $(-1,1)\times(-1,1).$ It is easy to write down the optimal transport map from \(U\) to \(V\):
\[
T(x,y)=
\begin{cases}
(x+1,y), & x>0,\\
(x-1,y), & x<0.
\end{cases}
\]
The corresponding Brenier potential is
\[
u(x,y)=\frac{x^2+y^2}{2}+|x|.
\]
Let \(\bar u\) denotes the minimal convex extension of \(u\). Because \(\bar u\) coincides with the supremum of the supporting planes of \(u\) over \(B_1\),
\[
\bar u(x,y)
=
|x|+
\sup_{a^2+b^2\leq 1}
\left\{
ax+by-\frac{a^2+b^2}{2}
\right\}.
\]
Setting \(r=\sqrt{x^2+y^2}\), this supremum evaluates to
\[
\begin{cases}
\dfrac{x^2+y^2}{2}, & r\leq 1,\\[6pt]
r-\dfrac12, & r\geq 1,
\end{cases}
\]
so that
\[
\bar u(x,y)
=
|x|+
\begin{cases}
\dfrac{x^2+y^2}{2}, & x^2+y^2\leq 1,\\[6pt]
\sqrt{x^2+y^2}-\dfrac12, & x^2+y^2\geq 1.
\end{cases}
\]
One can check that \(S_{\bar u}\) is the entire \(y\)-axis, hence $\Sigma=\{0\}\times\R$. For $|b|<1$, we have $\partial\bar u(0,b)=[-1,1]\times\{b\}$. Thus, by Legendre duality,
\[
\nabla v(a,b)=(0,b)
\qquad\text{for }(a,b)\in(-1,1)^2.
\]
Since the cavity $\operatorname{ch}(V)\setminus V$ coincides with $(-1,1)^2$ up to a null set, each horizontal segment of the cavity collapses onto a point of the singular diameter. Consequently,
\[
\mu=2\,\mathcal H^1\!\restriction_{\{0\}\times[-1,1]},
\qquad
\operatorname{spt}\mu=\{0\}\times[-1,1],
\]
where $\restriction$ denotes measure restriction: $\mu\restriction_A[B] := \mu[A \cap B]$ for all $B$. Thus the global equality \(\Sigma=\operatorname{spt}\mu\) breaks down, whereas Theorem~\ref{main} continues to hold:
\[
\Sigma\cap U=\operatorname{spt}\mu\cap U=\{0\}\times(-1,1).
\]
\end{ex}
Note that, within \(U\), the discontinuity set reduces to the diameter that separates the two transported halves. Beyond \(U\), on the other hand, the discontinuity set of the minimal convex extension persists along the entire \(y\)-axis, yet this exterior portion receives no phantom mass from \(\ch(V)\setminus V\). This is exactly the reason the identity has to be localized to the source: a global equality between \(\Sigma\) and \(\operatorname{spt}\mu\) is too much to hope for, and only the localized form \(\Sigma\cap U=\operatorname{spt}\mu\cap U\) survives.\\

In the next example,  we show existence of optimal transport map where $\operatorname{spt}\mu$ can entirely lie outside $U$ and it can even extend outside $\overline{U}$:
\begin{ex}\label{ex:explicit}
Define
\[
u(x,y)
=
\frac{(x+\frac32)^3}{6}
+
\frac{y^2}{2(x+\frac32)}
\qquad\text{in }U=B_1.
\]
Here $B_1$ is the unit ball in $\mathbb{R}^2$ centered at the origin.

\noindent A direct computation gives
\[
D^2u
=
\begin{pmatrix}
x+\frac32+\dfrac{y^2}{(x+\frac32)^3}
&
-\dfrac{y}{(x+\frac32)^2}
\\[1em]
-\dfrac{y}{(x+\frac32)^2}
&
\dfrac{1}{x+\frac32}
\end{pmatrix},
\qquad
\det D^2u=1.
\]
Since $x+\frac32>0$ on $\overline{B_1}$, for $e_2=(0,1)$ we have
\[
e_2^T D^2u\,e_2=\frac{1}{x+\frac32}>0.
\]
Moreover, $\det D^2u=1>0$, so the two eigenvalues of the symmetric matrix $D^2u$ have the same sign. They cannot both be negative by the preceding inequality. Hence both eigenvalues are positive, and therefore $D^2u>0$. Thus $u$ is strictly convex and smooth on a neighborhood of $\overline{B_1}$, and $\nabla u$ is a smooth diffeomorphism from $B_1$ onto \( P:=\nabla u(B_1). \) Since $\det D^2u=1$, we have \( (\nabla u)_\#(\chi_{B_1}\,dx)=\chi_P\,dy. \) As $\nabla u$ is the gradient of a convex function, it is the optimal transport map from $B_1$ to $P$ for the quadratic cost.\\

We now compute $P$ explicitly. Write
$$
(\xi,\eta)=\nabla u(x,y).
$$
Then
$$
\xi
=
\frac{(x+\frac32)^2}{2}
-
\frac{y^2}{2(x+\frac32)^2},
\qquad
\eta
=
\frac{y}{x+\frac32}.
$$
Set
$$
a=x+\frac32.
$$
Then $a>0$, and
$$
y=a\eta,
\qquad
\xi=\frac{a^2-\eta^2}{2}.
$$
Equivalently,
$$
a=\sqrt{\eta^2+2\xi}.
$$
The condition $(x,y)\in B_1$ becomes
$$
\left(a-\frac32\right)^2+a^2\eta^2<1.
$$
Therefore
$$
P
=
\left\{
(\xi,\eta):
\left(\sqrt{\eta^2+2\xi}-\frac32\right)^2
+
\eta^2(\eta^2+2\xi)
<1
\right\}.
$$

\begin{figure}[ht]
    \centering
    \includegraphics[width=0.8\textwidth]{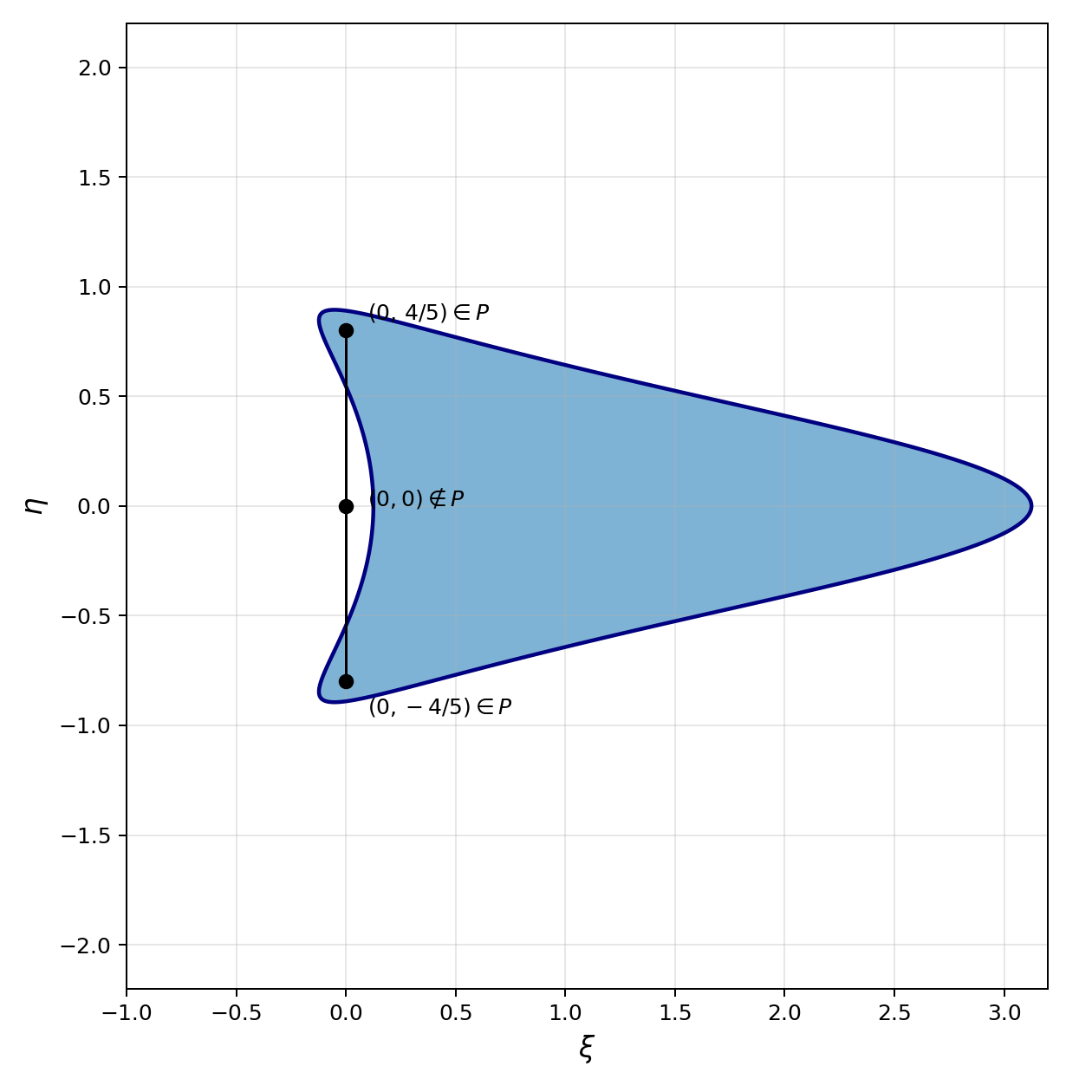}
    \caption{The nonconvex target region $P=\nabla u(B_1)$.}
    \label{fig:Vregion}
\end{figure}
\noindent We claim that $P$ is not convex. Indeed,
\[
\nabla u\left(-\frac{7}{10},\pm\frac{16}{25}\right)
=
\left(0,\pm\frac45\right),
\]
and both preimages belong to $B_1$. Hence
\[
\left(0,\pm\frac45\right)\in P.
\]
On the other hand, $(0,0)\notin P$: if $\nabla u(x,y)=(0,0)$, then $y=0$, and the first component of $\nabla u$ becomes
\[
\frac{(x+\frac32)^2}{2}>0,
\]
a contradiction. Therefore $P$ contains two points but not their midpoint, so $P$ is not convex. See Figure \ref{fig:Vregion} for an image of $P$.

\noindent Thus we have an explicit optimal transport map that is smooth in $U$, although the target $P$ is not convex. This means $\operatorname{spt}\mu$  is nonempty, but $\operatorname{spt}\mu\cap { \bar U}=\emptyset.$ {Since \(\mu(\partial B_1)=|\nabla { \bar u}(\partial B_1)|=0\), we also obtain
\[
\emptyset\neq
\spt\mu\cap(\mathbb R^2\setminus\overline{B_1})
\subset
\Sigma\cap(\mathbb R^2\setminus\overline{B_1}).
\]}
\end{ex}

\section{Global Partial Regularity}\label{partialsec}
First we give a proof of the global partial regularity for polygonal target: Theorem \ref{partial}.

\begin{proof}
{Let \(A:= S_{\bar u}\cap U\) be the singular set. Let \(E_i\) be the sides of \(V\), with \(E_i^\circ\) denoting the relative interior of \(E_i\). Define
\[
F_1:=\{x\in A:\partial\bar u(x)\text{ contains a vertex of }V\}. 
\]
The density bounds give \(\mathcal M_{\bar u}\ge c\,dx\) in \(U\) and \(\mathcal M_v\ge c\,dy\) in \(V\), where \(c=\lambda/\Lambda>0\). Thus \(\bar u\) and \(v\) are locally strictly convex in \(U\) and \(V\), respectively, by Alexandrov's planar 
differentiability theorem \cite{Al1}. 

The set \(F_1\) is finite. Indeed, otherwise, since \(V\) has finitely many vertices, there are distinct \(x,z\in F_1\) whose subgradients contain the same vertex. Then \(\bar u\) is affine on \([x,z]\). This contradicts local strict convexity near \(x\) inside \(U\), since a nontrivial initial portion of this segment lies in \(U\).

For \(x\in A\setminus F_1\), the extreme point argument in \cite{C1}, together with local strict convexity of \(v\) in \(V\), shows that every extreme point of \(\partial\bar u(x)\) lies on \(\partial V\). Moreover, each open side contains at most one such point. Indeed, if distinct \(p,q\in E_i^\circ\) belong to \(\partial\bar u(x)\), then \(v\) agrees with a supporting affine function \(\ell\) on \([p,q]\). Choose a small half-ball contained in \(V\), with diameter contained in \((p,q)\). Since \(v\ge\ell\) and \(v=\ell\) on the diameter, extending \(v\) by \(\ell\) to the opposite half-ball gives a convex function affine on the diameter. The lower bound \(\mathcal M_v\ge c\,dy\) on the  half-ball in $V$ then contradicts \cite[Lemma~6.1]{MR1} after rescaling. 

Since \(x\notin F_1\), the subdifferential contains no vertex. Being compact, convex, and not a singleton, it therefore has at least two extreme points lying in the relative interiors of distinct sides of \(V\).

For \(i<j\), set
\[
\Gamma_{ij}:=\{x\in A\setminus F_1:
\partial\bar u(x)\cap E_i^\circ\neq\emptyset,\
\partial\bar u(x)\cap E_j^\circ\neq\emptyset\}.
\]
Then
\[
A=F_1\cup\bigcup_{i<j}\Gamma_{ij}.
\]
The classes \(\Gamma_{ij}\) may overlap at points whose subgradients meet three or more open sides, {in which case \(\partial\bar u(x)\) will} not be a segment.

Fix \(i<j\) such that \(\Gamma_{ij}\neq\emptyset\). Since \(V\) is a polygon, two distinct sides are either disjoint or meet at a single endpoint. In either case, the directions of all the line segments joining \(E_i\) to \(E_j\) lie in a fixed cone whose opening is strictly smaller than \(\pi\). Thus, after rotating coordinates, vectors in this cone point uniformly in the vertical direction.

Let \(x,z\in\Gamma_{ij}\), and choose
\[
p_x,q_x\in\partial\bar u(x),\qquad
p_z,q_z\in\partial\bar u(z),
\]
where \(p_x,p_z\in E_i^\circ\) and \(q_x,q_z\in E_j^\circ\). By monotonicity of the subdifferential,
\[
(q_x-p_z)\cdot(x-z)\geq0,\qquad
(q_z-p_x)\cdot(x-z)\leq0.
\]
Both \(q_x-p_z\) and \(q_z-p_x\) belong to the fixed cone described above. Therefore the direction of \(x-z\), when \(x\neq z\), must stay a uniform positive angle away from the vertical direction.

It follows that \(\Gamma_{ij}\) is contained in the graph of a Lipschitz function, whose Lipschitz constant \(L_{ij}\) depends only on the two sides \(E_i\) and \(E_j\).  Extend this function with the same Lipschitz constant to the smallest closed interval containing its domain. The resulting graph is closed and therefore contains \(\overline{\Gamma_{ij}}\). Since \(U\) is bounded, this interval has length at most \(\operatorname{diam}(U)\). Hence
\[
\H^1\bigl(\overline{\Gamma_{ij}}\cap U\bigr)
\leq \sqrt{1+L_{ij}^2}\,\operatorname{diam}(U).
\]
It is clear that
\[
\Sigma\cap U
=
U\cap\overline A
\subset
F_1\cup\bigcup_{i<j}
\bigl(\overline{\Gamma_{ij}}\cap U\bigr).
\]
Since \(F_1\) is finite, \(\mathcal H^1(F_1)=0\). Summing the preceding estimates gives
\[
\mathcal H^1(\Sigma\cap U)
\le \sum_{i<j}\mathcal H^1(\overline{\Gamma}_{ij}\cap U)
\le C_V\,\operatorname{diam}(U),
\]
where \(C_V:=\sum_{i<j}\sqrt{1+L_{ij}^2}\) depends only on the geometry of \(V\). For empty classes, \(L_{ij}\) may be chosen from the same cone geometry.}
\end{proof}
\begin{rem}
{The finiteness of \(F_1\) was also established in \cite[Lemma~3.4]{CLL1}. We have included the elementary vertex-counting argument for completeness. The enlarged side-pair classes used above allow us to obtain the length estimate without classifying the multiple-contact points or proving their finiteness.}
\end{rem}

\begin{rem}\label{rem:explicit}
Let \(E_i=[a_i,b_i]\), \(i=1,\ldots,m\), be the sides of \(V\). For \(i<j\), let \(\theta_{ij}\in[0,\pi)\) denote the opening angle of the smallest closed convex cone containing all nonzero vectors \(q-p\), with \(p\in E_i\) and \(q\in E_j\). Then the preceding proof gives
\[
\H^1(\Sigma\cap U)\leq C_V\,\operatorname{diam}(U),
\]
where \(C_V\) can be chosen to satisfy
\[
C_V\leq\sum_{1\leq i<j\leq m}
\frac{1}{\cos(\theta_{ij}/2)}.
\]
Indeed, taking the bisector of each cone as the vertical direction in the monotonicity argument gives \(L_{ij}\leq\tan(\theta_{ij}/2)\). Thus the bound is explicitly determined by the vertices of \(V\).
\end{rem}

{Theorem~\ref{main} identifies the support of the phantom measure inside \(U\), while Theorem~\ref{partial} controls its length for polygonal targets. For uniform densities, we can further describe the measure itself.

Let \(A\) and \(F_1\) be as in the proof of Theorem~\ref{partial}, and set \(F:=F_1\cup F_2\), where
\[
F_2:=\{x\in A\setminus F_1:\partial\bar u(x)
\text{ meets three or more open sides of }V\}.
\]
The set \(F_1\) is finite by the  proof of Theorem \ref{partial}, and \(F_2\) is finite by \cite[Lemma~3.5]{CLL1}. Thus \(F\) is finite. For \(x\in A\setminus F\), the preceding proof shows that \(\partial\bar u(x)\) is a nontrivial segment joining the relative interiors of two distinct sides of \(V\). Moreover, \(A\setminus F\) is locally a smooth curve by \cite[Proposition~3.1]{CLL1}. We obtain the following description of the phantom measure, 
{in which the regularity of $\theta$ along these curves is related to the smoothness of the restriction of $u$ by Chiappori, McCann and Pass' theory of nested transport between unequal dimensions \cite{ChiapporiMcCannPass17}; (c.f. the appendix to \cite{McCannPass20}).}

\begin{cor}\label{cor:structure}
{For uniform densities and a polygonal target, we have
\[
\mu\restriction_U
=
\theta\,\mathcal H^1\restriction_{A\setminus F}
+
\sum_{x\in F}|\partial\bar u(x)|\,\delta_x,
\qquad
0\le\theta\in L^1(\mathcal H^1\restriction_{A\setminus F}),
\]
where \(A\) and \(F\) are defined above. }
\end{cor}

\begin{proof}
{Cover $A\setminus F$ by countably many compact regular arcs. On each such arc, take an injective arclength parametrization $\gamma:J\to A\setminus F$ and write
\[
\partial\bar u(\gamma(s))=[p(s),q(s)].
\]
By the proof of Proposition 3.1 of \cite{CLL1}, the endpoint maps $p,q$ are Lipschitz on $J$. Hence
\[
Y(s,t):=(1-t)p(s)+tq(s)
\]
is Lipschitz on $J\times[0,1]$.

Let $N\subset\gamma(J)$ be Borel with $\mathcal H^1(N)=0$, and set $E:=\gamma^{-1}(N)$. Since $\gamma$ is parametrized by arclength, $|E|=0$. Thus $E\times[0,1]$ has zero area, and its Lipschitz image satisfies
\[
|\partial\bar u(N)|
=
|Y(E\times[0,1])|=0.
\]
The identity \eqref{eq:phantom-subgradient} gives $\mu(N)=0$. The countable covering therefore yields
\[
\mu\restriction_{A\setminus F}
\ll\mathcal H^1\restriction_{A\setminus F}.
\]
By the Radon--Nikodym theorem, this restriction has a nonnegative density $\theta$, and
\[
\int_{A\setminus F}\theta\,d\mathcal H^1
=\mu(A\setminus F)
\le|\operatorname{ch}(V)\setminus V|<\infty.
\]
Finally, since \(A=S_{\bar u}\cap U\), Equation \eqref{eq:no-mass-on-closure-points} gives \(\mu(U\setminus A)=\mu(U\setminus S_{\bar u})=0\). Moreover, \(\mu(\{x\})=|\partial\bar u(x)|\) follows from the Monge--Amp\`ere identity. Combining these identities with the representation of \(\mu\) on \(A\setminus F\) completes the proof.}
\end{proof}

{As shown in Subsection \ref{sharp}, the exact equality in Theorem~\ref{main} can fail in higher dimensions. However, if the target \(V\) is ``nice enough'', one may still expect equality. The following remark establishes this for polyhedral targets in every dimension.}

\begin{rem}\label{rem:equality}
{For polyhedral targets, the exact support identity in Theorem~\ref{main} extends to every dimension:
\[
\Sigma\cap U=\spt\mu\cap U.
\]
Indeed, \cite[Theorem~3.1]{CLLLW} gives strict convexity of the Brenier potential in \(U\), without assuming convexity of \(U\). On every ball \(B\Subset U\setminus\spt\mu\), the phantom identity gives \(d\mathcal M\bar u=dx\). Caffarelli's interior regularity theory \cite{C1} therefore implies that \(\bar u\) is continuously differentiable in \(B\), so \(B\cap\Sigma=\emptyset\). The opposite inclusion follows from Theorem~\ref{main}.}
\end{rem}

\section{Future Directions}\label{future}
In this final section we list several open problems in this direction:
\begin{itemize}
    \item \textbf{Stability under domain and density perturbations:} How do the phantom measure and the geometry of its support change under small perturbations of the source and target domains and densities, particularly when uniform densities are replaced by nearby nonuniform densities?
    \item \textbf{Quantitative partial regularity for polytope target in higher dimensions:} Does Theorem~\ref{partial} extend to nonconvex polyhedral targets in dimensions \(n\geq 3\)? Under the same density bounds, can one establish
    \[
    \mathcal H^{n-1}(\Sigma\cap U)
    \le C_V\operatorname{diam}(U)^{n-1},
    \]
    with \(C_V\) depending only on the geometry of (V)? Such an estimate would complement the local regularity theorem of Chen, Li, Li, Liu and Wu \cite[Theorem~1.1]{CLLLW} by providing global measure control, including near the source boundary. More generally, which geometric assumptions on the target ensure that \(\mathcal H^{n-1}(\Sigma\cap U)<\infty\)?
    \item \textbf{Prove Theorem \ref{partial} for analytic $V$:} The idea in Section~\ref{partialsec} should extend to targets $V$ whose boundary consists of one connected analytic arc, along which the curvature vector points into the domain, together with finitely many line segments. In this setting, the geometry of the target should force the discontinuity set to consist of finitely many regular curves, each of finite length.  One can also expect to prove partial regularity in very symmetric setting as in Examples in  Rakshit and Sen \cite{RS1}. However, for a target with a general analytic boundary, it is not clear how to first prove that only finitely many regular curves occur. {For smooth source and target domains in $\R^n$, Chen, Li, Li, Liu and Wu \cite{CLLLW2} recently proved that the singular set has Hausdorff dimension at most \(n-1\). A further question is whether $\H^1(\Sigma\cap U)<\infty$ for targets with general analytic boundary in the plane.}
\end{itemize}

\end{document}